\documentclass[12pt]{article}
\usepackage{graphicx}
\usepackage{amsmath,amsthm,amssymb,enumerate}
\usepackage{euscript,mathrsfs}
\usepackage{color}
\usepackage{dsfont}
\usepackage[left=2cm,right=2cm,top=3.5cm,bottom=3.5cm]{geometry}
\usepackage{color}
\usepackage[framemethod=tikz]{mdframed}
\allowdisplaybreaks

\usepackage{soul}

\catcode`\@=11 \@addtoreset{equation}{section}

\catcode`\@=12

\allowdisplaybreaks

\newtheorem{Theorem}{Theorem}[section]
\newtheorem{Proposition}[Theorem]{Proposition}
\newtheorem{Lemma}[Theorem]{Lemma}
\newtheorem{Corollary}[Theorem]{Corollary}

\theoremstyle{definition}
\newtheorem{Definition}[Theorem]{Definition}

\newtheorem{Remark}[Theorem]{Remark}

\newcommand{\bTheorem}[1]{
\begin{Theorem} \label{T#1} }
\newcommand{\eT}{\end{Theorem}}

\newcommand{\bProposition}[1]{
\begin{Proposition} \label{P#1}}
\newcommand{\eP}{\end{Proposition}}

\newcommand{\bLemma}[1]{
\begin{Lemma} \label{L#1} }
\newcommand{\eL}{\end{Lemma}}

\newcommand{\bCorollary}[1]{
\begin{Corollary} \label{C#1} }
\newcommand{\eC}{\end{Corollary}}

\newcommand{\bRemark}[1]{
\begin{Remark} \label{R#1} }
\newcommand{\eR}{\end{Remark}}

\newcommand{\bDefinition}[1]{
\begin{Definition} \label{D#1} }
\newcommand{\eD}{\end{Definition}}

\newcommand{\Del}{\Delta_x}

\newcommand{\Ds}{\mathbb{D}_x}

\newcommand{\bfphi}{\boldsymbol{\varphi}}

\newcommand{\bFormula}[1]{
\begin{equation} \label{#1}}
\newcommand{\eF}{\end{equation}}

\newcommand{\vrn}{\vr_n}
\newcommand{\vun}{\vu_n} 
\newcommand{\vtn}{\vt_n} 
\newcommand{\chin}{\chi_n}

\newcommand{\Ov}[1]{\overline{#1}}

\newcommand{\DC}{C^\infty_c}

\newcommand{\aleq}{\stackrel{<}{\sim}}

\newcommand{\vr}{\varrho}

\newcommand{\vt}{\vartheta}
\newcommand{\vu}{\vc{u}}
\newcommand{\vm}{\vc{m}}

\newcommand{\vc}[1]{{\bf #1}}

\newcommand{\Div}{{\rm div}_x}
\newcommand{\Grad}{\nabla_x}

\newcommand{\dx}{\,{\rm d} {x}}

\newcommand{\dt}{\,{\rm d} t }

\newcommand{\intO}[1]{\int_{\Omega} #1 \ \dx}

\newcommand{\D}{{\rm d}}

\newcommand{\ep}{\varepsilon}

\newcommand{\br}{ \nonumber \\ }

\def\softd{{\leavevmode\setbox1=\hbox{d}%
          \hbox to 1.05\wd1{d\kern-0.4ex{\char039}\hss}}}
\definecolor{Cgrey}{rgb}{0.85,0.85,0.85}
\definecolor{Cblue}{rgb}{0.50,0.85,0.85}
\definecolor{Cred}{rgb}{1,0,0}
\definecolor{fancy}{rgb}{0.10,0.85,0.10}

\newcommand\Cbox[2]{%
    \newbox\contentbox%
    \newbox\bkgdbox%
    \setbox\contentbox\hbox to \hsize{%
        \vtop{
            \kern\columnsep
            \hbox to \hsize{%
                \kern\columnsep%
                \advance\hsize by -2\columnsep%
                \setlength{\textwidth}{\hsize}%
                \vbox{
                    \parskip=\baselineskip
                    \parindent=0bp
                    #2
                }%
                \kern\columnsep%
            }%
            \kern\columnsep%
        }%
    }%
    \setbox\bkgdbox\vbox{
        \color{#1}
        \hrule width  \wd\contentbox %
               height \ht\contentbox %
               depth  \dp\contentbox
        \color{black}
    }%
    \wd\bkgdbox=0bp%
    \vbox{\hbox to \hsize{\box\bkgdbox\box\contentbox}}%
    \vskip\baselineskip%
}

\mdfdefinestyle{MyFrame}{%
	linecolor=black,
	outerlinewidth=1pt,
	roundcorner=5pt,
	innertopmargin=-0.3\baselineskip,
	innerbottommargin=0.5\baselineskip,
	innerrightmargin=10pt,
	innerleftmargin=10pt,
	backgroundcolor=gray!20!white}

\date{}

\begin{document}


\title{A model of a non-isothermal two phase flow of compressible fluids\thanks{The authors thank the mobility project 8J20FR007 Barrande 2020 of collaboration between France and Czech Republic. The grantor in the Czech Republic is the Ministry of Education, Youth and Sports.}
}

\author{Eduard Feireisl$^{1}$\and M\u ad\u alina Petcu$^{2,3,4}$ \and Bangwei She$^{1,5}$  }


\maketitle

\medskip

\centerline{${}^1$Institute of Mathematics of the Academy of Sciences of the Czech Republic}
\centerline{\v Zitn\' a 25, CZ-115 67 Praha 1, Czech Republic}

\bigskip

\centerline{${}^2$Laboratoire de Math\' ematiques et Applications, UMR CNRS 7348 - SP2MI}

\centerline{Universit\' e de Poitiers, Boulevard Marie et Pierre Curie - T\' el\' eport 2}

\centerline{86962 Chasseneuil, Futuroscope Cedex,
	France}

\centerline{${}^3$The Institute of Mathematics of the Romanian Academy, Bucharest, Romania}

\centerline{and}

\centerline{${}^4$The Institute of Statistics and Applied Mathematics of the Romanian Academy, Bucharest, Romania}

\bigskip

\centerline{${}^5$Department of Mathematical Analysis, Charles University}
\centerline{Sokolovsk\' a 83, CZ-186 75 Praha 8, Czech Republic}

\begin{abstract}
	
We introduce a simple model of the time evolution of a binary mixture of compressible fluids including the thermal effects. Despite its apparent simplicity, the model is thermodynamically consistent admitting an entropy balance equation. We introduce a suitable weak formulation of the problem based on a combination of the entropy inequality with the total energy conservation principle. Finally, we show compactness of any bounded family of weak solutions and establish a global existence result.

\end{abstract}


\bigskip


\section{Model}
\label{M}

In \cite{FePePr} we have introduced a simple model of the motion of a mixture of two compressible fluids: 
\begin{subequations}\label{i1}
\begin{equation}\label{i1D}
\partial_t \vr + \Div (\vr \vu) = 0,  
\end{equation}
\begin{equation}\label{i1M}
\partial_t (\vr \vu) + \Div (\vr \vu\otimes \vu) = \Div \mathbb{T}(\vr, \chi, \Grad \chi, \Ds \vu),  
\end{equation}
\begin{equation}\label{i1C}
\partial_t \chi + \vu \cdot \Grad \chi =  \Del \chi - \frac{f(\chi) }{\partial \chi },	
\end{equation}
\end{subequations}
in a bounded domain $\Omega \subset R^d$, $d=2,3$, 
where $\vr$ is the fluid density, $\vu$ the bulk velocity, $\mathbb{T}$ the Cauchy stress, and $\chi$ an order parameter corresponding to the concentration difference.  Here, the Cauchy stress takes the form 
\begin{equation} \label{i2}
\mathbb{T} = \mathbb{S}(\Ds \vu) - \Big[ p(\vr) - f(\chi) \Big] \mathbb{I} - 
 \Big[ (\Grad \chi \otimes \Grad \chi) - \frac{1}{2} |\Grad \chi |^2 \mathbb{I} \Big],\ 
\Ds \vu = \frac{1}{2} \left( \Grad \vu + \Grad^t \vu \right), 
\end{equation}
with a Newtonian viscous stress 
\begin{equation} \label{i3}
	\mathbb{S} (\Ds \vu ) = \mu \left[  \Ds \vu - \frac{1}{d} \Div \vu \mathbb{I} \right] 
	+ 	 \eta \Div \vu \mathbb{I}.
	\end{equation}
The Allen--Cahn equation \eqref{i1C} depends only on the order parameter $\chi$ and the velocity $\vu$ as a consequence of our choice of the free energy, namely
\[
E_{\rm free} = \frac{1}{2} |\Grad \chi|^2 + \vr e(\vr) + f(\chi), 
\]
with $e$ independent of $\chi$.

The goal of the present paper is to incorporate the effect of temperature changes  into the model following the same philosophy of separation 
of the thermodynamic variables from the order parameter in the pressure -- internal energy equation of state. To this end, we suppose that both 
the pressure $p$ and the viscosity coefficients $\mu$ and $\eta$ depend on the (absolute) temperature $\vt$. In addition, we also assume the viscosity may be different for each mixture component, meaning $\mu = \mu(\vt, \chi)$, 
$\eta = \eta(\vt, \chi)$. Motivated by \eqref{i1C}, we suppose the effective pressure to be given as 
\[
p_{\rm eff}(\vr, \vt, \chi) = p(\vr, \vt) - f(\chi).
\] 

\subsection{Energy balance}

Our goal is  to show that the model is thermodynamically consistent. First,  we write the momentum equation \eqref{i1M} in the form 
\begin{equation} \label{M2}
	\partial_t (\vr \vu) + \Div (\vr \vu \otimes \vu) + \Grad p(\vr, \vt) = \Div \mathbb{S} - 
	\Div \left( \Grad \chi \otimes \Grad \chi - \frac{1}{2} |\Grad \chi|^2 \mathbb{I} \right) + \Grad f(\chi),
\end{equation}
with 
\begin{equation} \label{M3}
	\mathbb{S} = \mu (\vt, \chi) \left(  \Ds \vu - \frac{1}{d} \Div \vu \mathbb{I} \right) + \eta (\vt, \chi) \Div \vu \mathbb{I}. 
\end{equation}

Next, multiplying \eqref{M2} with $\vu$ we get 
\begin{align}
\partial_t \left[ \frac{1}{2} \vr |\vu|^2 \right] &+ \Div \left( \left[ \frac{1}{2} \vr |\vu|^2 + p - f \right] \vu \right) 
- \Div \left(\mathbb{S} \cdot \vu \right) \br &= (p - f) \Div \vu - \mathbb{S} : \Ds \vu - \Div \left( \Grad \chi \otimes \Grad \chi - \frac{1}{2} |\Grad \chi|^2 \mathbb{I} \right) \cdot \vu.
\label{M4}
\end{align}
Seeing that 
\[
\Div \left( \Grad \chi \otimes \Grad \chi - \frac{1}{2} |\Grad \chi|^2 \mathbb{I} \right) = 
\Del \chi \Grad \chi
\]
we may consider an equation for the order parameter $\chi$ of Allen--Cahn type: 
\begin{equation} \label{M5}
	\partial_t \chi + \vu \cdot \Grad \chi = \Del \chi - \frac{\partial f}{\partial \chi }
	\end{equation}

Multiplying \eqref{M5} on $\Del \chi$ and adding the result to \eqref{M4} yields 
\begin{align}
	\partial_t \left[ \frac{1}{2} \vr |\vu|^2 + \frac{1}{2} |\Grad \chi|^2 \right] &+ \Div \left( \left[ \frac{1}{2} \vr |\vu|^2 + ( p - f) \right] \vu \right) - \Div (\partial_t \chi \Grad \chi ) 
	- \Div \left(\mathbb{S} \cdot \vu \right) \br &= ( p - f) \Div \vu - \mathbb{S} : \Ds \vu - \Del \chi \left( \Del \chi - \frac{\partial f}{\partial \chi } \right).
	\label{M6}
\end{align}
Similarly, multiplying \eqref{M5} on $\frac{\partial f}{\partial \chi}$ we deduce 
\begin{multline}	\label{M7}
	\partial_t \left[ \frac{1}{2} \vr |\vu|^2 + \frac{1}{2} |\Grad \chi|^2 \right] + \Div \left( \left[ \frac{1}{2} \vr |\vu|^2 + (p - f) \right] \vu \right) - \Div (\partial_t \chi \Grad \chi ) 
	- \Div \left(\mathbb{S} \cdot \vu \right) 
\\	+ \frac{\partial f}{\partial \chi} \partial_t \chi +  \vu \cdot \Grad \chi \frac{\partial f}{\partial \chi } 
	 = (p  - f ) \Div \vu - \mathbb{S} : \Ds \vu -  \left( \Del \chi -  \frac{\partial f}{\partial \chi } \right)^2.
\end{multline}
Consequently, the mechanical energy balance reads
\begin{multline}	\label{M9}
	\partial_t \left[ \frac{1}{2} \vr |\vu|^2 + \frac{1}{2} |\Grad \chi|^2 + f(\chi) \right] + \Div \left( \left[ \frac{1}{2} \vr |\vu|^2 + p \right] \vu \right) - \Div (\partial_t \chi \Grad \chi ) - \Div \left(\mathbb{S} \cdot \vu \right)    
\\	=  p \Div \vu - \mathbb{S} : \Ds \vu -  \left( \Del \chi -  \frac{\partial f(\chi) }{\partial \chi } \right)^2.
\end{multline}

Introducing the internal energy $e = e(\vr, \vt)$ we write the First Law of Thermodynamics in the form 
\begin{align} 	\label{M12}
	\partial_t (\vr e) +\Div (\vr e \vu) + \Div \vc{q} = 
	- p(\vr, \vt) \Div \vu + \mathbb{S} : \Ds \vu + \left( \Del \chi -  \frac{\partial f(\chi) }{\partial \chi } \right)^2,
\end{align}
where $\vc{q}$ is the internal energy flux. Summing up \eqref{M9}, \eqref{M12} we therefore obtain the \emph{total energy balance}
\begin{multline}	\label{M13}
	\partial_t \left[ \frac{1}{2} \vr |\vu|^2 + \frac{1}{2} |\Grad \chi|^2 + \vr e + f(\chi) \right] + \Div \left( \left[ \frac{1}{2} \vr |\vu|^2 + \vr e + p \right] \vu \right) 
\\	- \Div (\partial_t \chi \Grad \chi )  - \Div \left(\mathbb{S} \cdot \vu \right)  + 	\Div \vc{q} = 0.  
\end{multline}

\subsection{Entropy}

The total energy balance \eqref{M13} is not suitable for a weak formulation as the convective fluxes therein are not controlled by the available {\it a priori} bounds. Therefore it is more convenient to introduce the entropy $s$ through 
\emph{Gibbs' equation} 
\begin{equation} \label{M11}
	\vt D s = D e + p D \left( \frac{1}{\vr} \right).
\end{equation}

Dividing the internal energy balance \eqref{M12} on $\vt$ and performing a routine manipulation we obtain the \emph{entropy equation} 
\begin{align} 
	\partial_t (\vr s) +\Div (\vr s \vu) + \Div \left( \frac{ \vc{q} }{\vt } \right) = 
	\frac{1}{\vt} \left[ \mathbb{S} : \Ds \vu + \left( \Del \chi -  \frac{\partial f(\chi) }{\partial \chi } \right)^2 - \frac{\vc{q} \cdot \Grad \vt }{\vt} 
	\right].
	\label{M14}
\end{align}
The right hand side of \eqref{M14} is the entropy production rate, and, in accordance with the Second Law of Thermodynamics, it is non--negative for any physically admissible process.

\subsection{Boundary conditions}

For the sake of simplicity, we  consider the \emph{no--slip} boundary conditions for the velocity, 
\begin{equation} \label{M17}
	\vu|_{\partial \Omega} = 0,
	\end{equation}
together with the Neumann boundary conditions for the order parameter and the heat flux, 
\begin{equation} \label{M18}
	\Grad \chi \cdot \vc{n}|_{\partial \Omega} = 0,\ 
	\vc{q} \cdot \vc{n}|_{\partial \Omega} = 0.
	\end{equation}
More general inhomogeneous boundary conditions can be handled by the approach developed in \cite{ChauFei} and 
\cite{FeiNov20}. 

Accordingly, the system is closed (mechanically and thermally insulated) conserving the total mass, 
\begin{equation} \label{M19}
	\frac{\D }{\dt} \intO{ \vr } = 0, 
	\end{equation}
and energy 
\begin{equation} \label{M20}
		\frac{\D }{\dt} \intO{ \left[ \frac{1}{2} \vr |\vu|^2 + \frac{1}{2} |\Grad \chi|^2 + \vr e + f(\chi) \right] } = 0,
	\end{equation}	
while the total entropy is non--decreasing, 
\[
	\frac{\D }{\dt} \intO{ \vr s } \geq 0.
\]

\subsection{Weak formulation}

Following the general approach developed in \cite{FeNo6A}, we relax the entropy equation
to \emph{inequality} and augmented the system by the (integrated) total energy balance. 
Accordingly, we consider the following problem for the unknowns $(\vr, \vt, \chi, \vu)$:

\begin{mdframed}[style=MyFrame]
\begin{equation}	\label{M15}
\begin{aligned} 
	&\partial_t \vr  + \Div (\vr \vu) = 0, \\
\partial_t (\vr \vu) + \Div (\vr \vu \otimes \vu) + \Grad p & = \Div \mathbb{S} - 
\Div \left( \Grad \chi \otimes \Grad \chi - \frac{1}{2} |\Grad \chi|^2 \mathbb{I} \right) + \Grad f,	
  \\
\partial_t \chi + \vu& \cdot \Grad \chi = \Del \chi - \frac{\partial f(\chi)}{\partial \chi }, 
 \\
\partial_t (\vr s) +\Div (\vr s \vu) + \Div \left( \frac{ \vc{q} }{\vt } \right) &\geq 
\frac{1}{\vt} \left[ \mathbb{S} : \Ds \vu + \left( \Del \chi -  \frac{\partial f(\chi) }{\partial \chi } \right)^2 - \frac{\vc{q} \cdot \Grad \vt }{\vt} 
\right],  
 \\
\frac{\D }{\dt } \intO{ \bigg[ \frac{1}{2} \vr |\vu|^2 &+ \frac{1}{2} |\Grad \chi|^2 + \vr e + f(\chi) \bigg] } 
 = 0.
	\end{aligned}
\end{equation}

\end{mdframed}

The above problem is suitable for a weak formulation and \emph{compatible} with the original system of equations as soon as all quantities 
are sufficiently smooth. By ``compatible" we mean that any sufficiently smooth $(\vr, \vt, \chi, \vu)$ satisfy the internal energy equation \eqref{M12}, and, consequently, the entropy equation \eqref{M14} as well as the total energy balance 
\eqref{M13}. Indeed, this can be shown exactly as in  \cite[Chapter 2]{FeNo6A}.

Finally, let us point out that
models of mixtures are mostly studied in the ``incompressible'' framework, meaning with the equation of continuity 
replaced by the solenoidality condition $\Div \vu = 0$ and the pressure (gradient) term in the momentum equation being implicitly determined by the motion, see e.g. Abels et al. \cite{AbeGarWeb,AbeTer}, the survey chapter by 
Abels and Garcke \cite{AbeGar}, and the references cited therein. A thermodynamically consistent model in the fully compressible setting was proposed by Blesgen \cite{Blesgen} and its ``barotropic'' version  was analysed in \cite{FePeRoSc}. We refer to 
Chen and Guo \cite{ChenGuo}, 
Kotschote~\cite{Kot1,Kot5,Kot4,Kot3}, and  Zhao et al.~\cite{ZGH} for other available analytical results.

\section{Main result}
\label{MR}

In this section, we introduce the principal hypotheses concerning the constitutive relations, define the concept of weak solution, and state our main result. 

\subsection{Hypotheses}

The hypotheses imposed on the pressure equation of state are inspired by \cite{FeNo6A}. We consider 
\begin{equation} \label{MR1}
	p(\vr, \vt) = \vt^{5/2} P \left( \frac{ \vr }{\vt^{3/2}} \right) + \frac{a}{3} \vt^4,\ a > 0,
\end{equation}
with $P \in C^1[0,\infty)$, where the term $\frac{a}{3} \vt^4$ represents the effect of thermal radiation. The compactification effect of this term necessary for the existence theory is largely discussed in \cite[Chapter 2,3]{FeNo6A}.
In accordance with Gibbs' equation (\ref{M11}), we get
\begin{align} 
	e(\vr, \vt) &= \frac{3}{2} \frac{\vt^{5/2}}{\vr} P \left( \frac{ \vr }{\vt^{3/2}} \right) + \frac{a}{\vr} \vt^4,
\br 
	s(\vr, \vt) &= \mathcal{S} \left( \frac{\vr} {\vt^{3/2}} \right) + \frac{4a}{3} \frac{\vt^3}{\vr},
\br
	\mathcal{S}'(Z) &= - \frac{3}{2} \frac{ \frac{5}{3} P(Z) - P'(Z) Z }{Z^2}.
	\label{MR2}
\end{align}
In addition, we assume
\begin{equation} \label{MR3}
	P'(Z) > 0 \ \mbox{for any}\ Z \geq 0,\ \frac{ \frac{5}{3} P(Z) - P'(Z) Z }{Z} > 0 \ \mbox{for any}\ Z > 0.
\end{equation}
In particular,  the function $Z \mapsto P(Z) / Z^{5/3}$ is decreasing, and we suppose
\begin{equation} \label{MR4}
	\lim_{Z \to \infty} \frac{ P(Z) }{Z^{5/3}} = p_\infty > 0.
\end{equation}
Note that the hypotheses \eqref{MR3}, \eqref{MR4} reflect the \emph{thermodynamic stability} - positive compressibility and positive specific heat at constant volume.
Finally, we suppose a technical but physically grounded hypothesis (see \cite[Chapter 2]{FeNo6A})
\begin{equation} \label{MR5}
	P(0) = 0,\ \frac{ \frac{5}{3} P(Z) - P'(Z) Z }{Z} < c \ \mbox{for all}\ Z > 0.
\end{equation}

The viscous stress $\mathbb{S}$ is given by Newton's rheological  law
\eqref{i3}, where the viscosity coefficients $\mu$ and $\eta$ are continuously differentiable functions of $\vt$, 
$\chi$ satisfying 
\begin{equation} \label{MR6}
	\underline{\mu} (1 + \vt) \leq \mu(\vt, \chi) \leq \Ov{\mu} (1 + \vt),\ | \nabla_{\vt, \chi} \mu(\vt, \chi) | < c,
\end{equation}
\begin{equation} \label{MR7}
	0 \leq \eta (\vt, \chi) \leq \Ov{\eta} (1 + \vt).
\end{equation}

Similarly, the heat flux is given by \emph{Fourier's law}, 
\begin{equation} \label{MR8}
	\vc{q} = - \kappa \Grad \vt,
	\end{equation}
where $\kappa = \kappa (\vt, \chi)$ is a differentiable function of the temperature, 
\begin{equation} \label{MR9}
	\underline{\kappa} (1 + \vt^3) \leq \kappa (\vt, \chi) \leq \Ov{\kappa}(1 + \vt^3).
\end{equation}	

The potential $f$ in the Allen--Cahn equation \eqref{M5} is chosen in such a way that the ``pure'' states 
$\chi = \pm 1$ are stationary solutions, 
\begin{equation} \label{MR10}
\partial_{\chi}	f(\pm 1) = 0.
\end{equation}
In particular, applying the standard parabolic maximum principle to \eqref{M5} we may infer that the order parameter 
remains confined to its natural range 
\begin{equation} \label{MR11}
	-1 \leq \chi(t,x) \leq 1 \ \mbox{as long as}\ -1 \leq \inf \chi(0,x) \leq \sup \chi(0,x) \leq 1.
	\end{equation}

\subsection{Weak solutions}

We suppose that 
\[
\Omega \subset R^d,\ d=2,3 \ \mbox{is a bounded domain.}
\]

\begin{mdframed}[style=MyFrame]
	\begin{Definition}[Weak solution] \label{DMR1}
		We say that $(\vr, \vt, \chi, \vu)$ is a \emph{weak solution} of the Navier--Stokes--Fourier--Allen--Cahn (NSFAC) system \eqref{M15}, with the boundary conditions \eqref{M17}, \eqref{M18}, and the initial conditions 
		\[
		\vr(0, \cdot) = \vr_0,\ \vr \vu (0, \cdot) = (\vr \vu)_0,\ \vr s (\vr, \vt)(0, \cdot) = \vr_0 s(\vr_0, \vt_0),\ 
		\chi(0) = \chi_0
		\] 
		if the following hold:
		
		\begin{itemize}
			\item {\bf Integrability, regularity.} 
			\[
			\vr \in C_{\rm weak}([0,T]; L^{\frac{5}{3}}(\Omega)),\ \vr \geq 0, 
			\]
			\[
			\vu \in L^2(0,T; W^{1,2}_0 (\Omega; R^d)),\ \vr \vu \in C_{\rm weak}(0,T; L^{\frac{5}{4}}(\Omega; R^d)),
			\]
			\[
			\vt \in L^\infty(0,T; L^4(\Omega)),\ \vt^{\frac{3}{2}} \in L^2(0,T; W^{1,2}(\Omega)), 
			\]
			\[
			\vt > 0 \ \mbox{a.a. in}\ (0,T) \times \Omega,\ 
			\log(\vt) \in L^2(0,T; W^{1,2}(\Omega)),
			\]
			\[
			\chi \in L^2(0,T; W^{2,2}(\Omega)),\ \partial_t \chi \in L^2(0,T; L^2(\Omega)), 
			\]
			\[
			-1 \leq \chi(t,x) \leq 1 \ \mbox{for a.a.}\ (t,x) \in (0,T) \times \Omega.
			\]
			
			\item {\bf Equation of continuity.}
			
			The integral identity 
			\begin{equation} \label{MR12}
			\int_0^T \intO{ \Big[ \vr \partial_t \varphi + \vr \vu \cdot \Grad \varphi \Big] } = - \intO{ \vr_0 
				\varphi (0, \cdot) }
			\end{equation}
			holds for any $\varphi \in C^1_c([0,T) \times \Ov{\Omega})$.
			
			\item {\bf Momentum equation.} 
			
			The integral identity 
			\begin{align} 
				\int_0^T \intO{ \Big[ &\vr \vu \cdot \partial_t \bfphi + \vr \vu \otimes \vu : \Grad \bfphi + p \Div \bfphi \Big] } \dt = \int_0^T \intO{ \mathbb{S} : \Grad \bfphi } \dt \br 
				&- \int_0^T \intO{ \left[ \Grad \chi \otimes \Grad \chi - \frac{1}{2} |\Grad \chi |^2 \mathbb{I} \right] : \Grad \bfphi } \dt  + \int_0^T \intO{ f \Div \bfphi } \dt \br &= - 
				\intO{ (\vr \vu)_0 \cdot \bfphi (0, \cdot) }
				\label{MR13}
			\end{align}
		holds for any $\bfphi \in C^1_c([0,T) \times \Omega; R^d)$.

		\item {\bf Allen--Cahn equation for order parameter.}
		\begin{align} 
			\partial_t \chi + \vu \cdot \Grad \chi &= \Del \chi - \frac{\partial f}{\partial \chi } 
			\ \mbox{a.a. in}\ (0,T) \times \Omega,\br \Grad \chi \cdot \vc{n}|_{\partial \Omega} &= 0,\  \chi(0, \cdot) = \chi_0.
			\label{MR14}
			\end{align}
		
		\item {\bf Entropy inequality.}
		
		The inequality 
		\begin{align} 
		- \int_0^T &\intO{ \Big[ \vr s \partial_t \varphi +\vr s \vu \cdot \Grad \varphi  + \left( \frac{ \vc{q} }{\vt } \right) \cdot \Grad \varphi \Big] } \dt \br &\geq 
		\int_0^T \intO{ \frac{\varphi}{\vt} \left[ \mathbb{S} : \Ds \vu + \left( \Del \chi -  \frac{\partial f(\chi) }{\partial \chi } \right)^2 - \frac{\vc{q} \cdot \Grad \vt }{\vt} 
		\right]	} \dt \br &+ \intO{ \vr_0 s(\vr_0, \vt_0) \varphi (0, \cdot) } \label{MR15}
			\end{align}
		holds for any $\varphi \in C^1_c ([0,T) \times \Ov{\Omega})$, $\varphi \geq 0$.

		\item {\bf Total energy conservation.}
		
		\begin{align} 
			& \intO{ \left[ \frac{1}{2} \vr |\vu|^2 + \frac{1}{2} |\Grad \chi |^2 + \vr e + f(\chi) \right] (\tau, \cdot) } \br 
			&= \intO{ \left[ \frac{1}{2} \frac{|(\vr \vu)_0|^2}{\vr_0} + \frac{1}{2} |\Grad \chi_0 |^2 + \vr_0 e(\vr_0, \vt_0) + f(\chi_0) \right] }
			\label{MR16}
			\end{align}
		for a.a. $0 \leq \tau \leq T$.
		
	\end{itemize}
		
		\end{Definition} 
	
	\end{mdframed}

\subsection{Main result}
	
	Finally, we are ready to formulate our main result.
	
	\begin{mdframed}[style=MyFrame]
		
		\begin{Theorem}[Global existence] \hfill \label{TMR1}

			Let $\Omega \subset R^d$, $d=2,3$ be a bounded domain of class $C^{2 + \nu}$. 
			Let the thermodynamic functions $p$, $e$, $s$ be given by \eqref{MR1}--\eqref{MR5}, let 
			the transport coefficients $\mu$, $\eta$, $\kappa$ satisfy \eqref{MR6}, \eqref{MR7}, \eqref{MR9}, and let 
			$f$ satisfy \eqref{MR10}.  
			Let the initial data 
			$(\vr_0, \vt_0, \chi_0, (\vr \vu)_0)$ belong to the class  
			\begin{align} 
				\vr_0  &\in L^{\frac{5}{3}} (\Omega),\ \vr_0 \geq 0,\ \frac{(\vr \vu)_0 }{\vr_0} \in L^1(\Omega), \br 
				\vt_0 & > 0 \ \mbox{a.a. in}\ \Omega,\ \vr_0 e(\vr_0, \vt_0),\ \vr_0 s(\vr_0 , \vt_0) \in L^1(\Omega), 
				\br 
				-1 &\leq \chi_0 \leq 1,\ \chi_0 \in W^{1,2}(\Omega).
				\label{MR17}  
				\end{align}
			
			\medskip
			
			Then for any $T > 0$, the NSFAC system admits a weak solution $(\vr, \vt, \chi, \vu)$ in $(0,T) \times \Omega$ in the sense of Definition \ref{DMR1}.

			\end{Theorem}

		\end{mdframed} 
	
	\begin{Remark}[Domain regularity] \label{RMR1}
		
		The hypothesis on regularity of the spatial domain can possibly be relaxed. Note, however, the full elliptic estimates 
		\[
		\Del \chi \in L^2,\ \Grad \chi \cdot \vc{n}|_{\partial \Omega} = 0 \ \Rightarrow \ 
		\chi \in W^{2,2}(\Omega)
		\]
		are needed to recover regularity of the order parameter.
		
		\end{Remark}
	
	The rest of the paper is devoted to the proof of Theorem \ref{TMR1}. In view of the existing theory for the Navier--Stokes--Fourier system \cite[Chapter 3]{FeNo6A}, we focus on {\it a priori} bounds and compactness of bounded sets of weak solutions (weak sequential stability). Finally, we introduce a family of approximate systems, similar to \cite{FeNo6A} to construct the weak solution, the existence of which is claimed in Theorem \ref{TMR1}.
	
	\section{A priori bounds}
	
	\label{ap}
	
	{\it A priori} bounds available for the NSFAC system are, to certain extent, similar to those for the Navier--Stokes--Fourier system. 
	
	\subsection{Boundedness of the order paremater} 
	
	As already pointed out, the standard parabolic maximum principle applied to the Allen--Cahn system yields 
	\begin{equation} \label{ap1}
-1 \leq \chi (t,x) \leq 1 \ \mbox{for}\ t \in (0,T),\ x \in \Omega.
\end{equation}

\subsection{Energy/entropy estimates} 

Similarly to \cite{ChauFei}, we introduce the \emph{ballistic energy} functional 
\[
E_{\Ov{\vt}} (\vr, \vt, \chi, \vu) = \frac{1}{2} \vr |\vu|^2 + \frac{1}{2} |\Grad \chi |^2 + \vr e + f(\chi) - 
\Ov{\vt} \vr s, 
\]
where $\Ov{\vt} > 0$ is a positive constant. Multiplying the entropy inequality \eqref{MR15} by $\Ov{\vt}$, integrating over $\Omega$, and adding the result to the total energy balance \eqref{MR16}, we obtain 
\begin{multline}	\label{ap2}
	\intO{ E_{\Ov{\vt}} (\vr, \vt, \chi, \vu) (\tau, \cdot) } + 
	\int_0^\tau \intO{ \frac{\Ov{\vt}}{\vt} \left[ \mathbb{S} : \Ds \vu + \left( \Del \chi -  \frac{\partial f(\chi) }{\partial \chi } \right)^2 - \frac{\vc{q} \cdot \Grad \vt }{\vt} \right] } \dt  
	\\ 
	\leq \intO{ \left[ \frac{1}{2} \frac{|(\vr \vu)_0|^2}{\vr_0} + \frac{1}{2} |\Grad \chi_0 |^2 + \vr_0 e(\vr_0, \vt_0) + f(\chi_0) - \Ov{\vt} \vr_0 s(\vr_0, \vt_0) \right] }
	\end{multline}
for $0 \leq \tau \leq T$. 

In view of \eqref{ap1}, the function $f(\chi)$ is bounded and we may deduce the following estimates: 
\begin{align} 
	\sup_{\tau \in [0,T]} \intO{ \left[ \frac{1}{2} \vr |\vu|^2 + \frac{1}{2} |\Grad \chi |^2 + \vr e + 
		\vr |s(\vr, \vt)| \right] } &\aleq 1, \br 
	\int_0^T 
	\intO{ \left[ \frac{\mu(\vt, \chi)}{\vt} \left|\Ds \vu - \frac{1}{d} \Div \vu \mathbb{I} \right|^2 + 
		\frac{\kappa (\vt, \chi) |\Grad \vt |^2 }{\vt^2} + \frac{1}{\vt} \left( \Del \chi -  \frac{\partial f(\chi) }{\partial \chi } \right)^2 \right] } \dt &\aleq 1. 
	\label{ap3}
	\end{align}
Using the structural restrictions imposed in the hypotheses \eqref{MR1}--\eqref{MR9}, the above estimates give rise to the uniform bounds:
\begin{align} 
	\sup_{\tau \in [0,T]} \| \vr (\tau, \cdot) \|_{L^{\frac{5}{3}}(\Omega) } &\aleq 1, \br 
		\sup_{\tau \in [0,T]} \| \vt (\tau, \cdot) \|_{L^{4}(\Omega) } &\aleq 1 , \br
\sup_{\tau \in [0,T]} \| \vr \vu (\tau, \cdot) \|_{L^{\frac{5}{4}}(\Omega; R^d) } &\aleq 1, \br
\sup_{\tau \in [0,T]} \| \chi (\tau, \cdot) \|_{W^{1,2}(\Omega) } &\aleq 1, \br
\int_0^T \| \vu \|_{W^{1,2}_0 (\Omega; R^d) }^2 \dt &\aleq 1, \br
\int_0^T \| \vt^{\frac{3}{2}} \|^2_{W^{1,2}(\Omega) } \dt &\aleq 1, \br 
\int_0^T \| \Grad \log \vt \|^2_{L^2(\Omega; R^d)} &\aleq 1.	
	\label{ap4}
	\end{align}

\subsection{Gagliardo--Nirenberg inequality and the Allen--Cahn equation} 

As a matter of fact, the inequality \eqref{ap3} yields certain bounds on 
$\Del \chi$ depending on integrability of $\vt$. However, a better bound can be deduced directly from the 
Allen--Cahn equation. As $\chi$ satisfies the Neumann boundary condition, we can use the Gagliardo--Nirenberg interpolation inequality, 
\begin{equation} \label{ap5}
	\| \Grad \chi \|_{L^4(\Omega; R^d)}^2 \leq \| \chi \|_{L^\infty (\Omega) } \| \Del \chi \|_{L^2(\Omega)}.
	\end{equation}

Next, multiplying the Allen--Cahn equation \eqref{MR14} on $\Del \chi$ and integrating by parts give rise to
\begin{equation} \label{ap6}
\frac{ \D }{\dt} \intO{ \frac{1}{2} |\Grad \chi |^2 } + \intO{ |\Del \chi |^2 } = \intO{ \vu \cdot \Grad \chi \Del \chi } + \intO{ \frac{\partial f (\chi) }{\partial \chi } \Del \chi }. 
\end{equation}
Now, we rewrite the first term on the right hand side as
\begin{align}
\intO{ \vu \cdot \Grad \chi \Del \chi } &= \intO{ \vu \cdot \Div \left( \Grad \chi \otimes \Grad \chi - \frac{1}{2} 
	|\Grad \chi |^2 \mathbb{I}  \right) } \br 
&= - \intO{ \Grad \vu : \left( \Grad \chi \otimes \Grad \chi - \frac{1}{2} 
	|\Grad \chi |^2 \mathbb{I}  \right) }. 
\label{ap5a}
\end{align}
Consequently, by virtue of  \eqref{ap1}, \eqref{ap5}, \eqref{ap5a}, the uniform bounds \eqref{ap4}, and H\" older's inequality,
\[
	\left| \intO{ \vu \cdot \Grad \chi \Del \chi } \right|  \aleq \| \vu \|_{W^{1,2}_0 (\Omega; R^d)} 
	\| \Grad \chi \|_{L^4(\Omega; R^d)}^2 
\aleq \| \vu \|_{W^{1,2}_0 (\Omega; R^d)} \| \Del \chi \|_{L^2(\Omega)}.
\]
Thus going back to \eqref{ap6} we may infer that 
\[
\int_0^T \intO{ |\Del \chi |^2 } \dt \aleq 1;
\]
whence, by standard elliptic estimates, 
\begin{equation} \label{ap7}
	\int_0^T \| \chi \|^2_{W^{2,2}(\Omega) } \dt \aleq 1.
	\end{equation}

\subsection{Pressure estimates}

In order to guarantee that all terms in the integral formulation \eqref{MR12}--\eqref{MR16} are bounded in a reflexive space $L^r$, $r > 1$, we need to control the pressure $p(\vr, \vt)$, more specifically, its ``elastic'' component  that is proportional to $\vr^{\frac{5}{3}}$. This step is well understood in the context of the Navier--Stokes--Fourier system (see e.g. \cite[Chapter 3, Section 3.6.3]{FeNo6A}) and easy to extend on the present setting given the uniform bounds \eqref{ap1} and \eqref{ap4}. We will come to this issue later in the forthcoming Section \ref{pe}.

\section{Compactness}
\label{c}

Our goal is to show compactness (weak sequential stability) property for any family of weak solutions 
$\{ \vrn, \vtn, \chin, \vun \}_{n \geq 0}$ of the NSFAC system satisfying the uniform bounds established in 
Section \ref{ap}, meaning emanating from bounded energy/entropy initial data.
In addition to Definition \ref{DMR1}, we also assume that the equation of continuity \eqref{MR12} is satisfied in the renormalized sense: 
\begin{align}
	\int_0^T &\intO{ \left[ b(\vr_n) \partial_t \varphi + b(\vr_n) \vun \cdot  \Grad \varphi + 
		\Big( b\vrn) - b'(\vrn) \vrn \Big) \Div \vun \varphi \right] } \dt \br &= 
	- \intO{ b(\vr_{n}(0, \cdot)) \varphi (0, \cdot) } 
	\label{renorm}
	\end{align}
for any $\varphi \in C^1_c ([0,T) \times \Ov{\Omega})$ and any $b \in C^1(R)$, $b' \in C_c(R)$.

In view of the uniform bounds \eqref{ap4}, we may extract a suitable subsequence such that 
\begin{align}
	\vrn &\to \vr \ \mbox{in}\ C_{\rm weak}([0,T]; L^{\frac{5}{3}}(\Omega)), \br
\vt_n &\to \vt \ \mbox{weakly-(*) in}\ L^\infty(0,T; L^4(\Omega)), \br
\chin &\to \chi \ \mbox{in, say,}\ L^2((0,T) \times \Omega),\ -1 \leq \chi \leq 1, \br
\vun &\to \vu \ \mbox{weakly in}\ L^2(0,T; W^{1,2}_0 (\Omega; R^d)).
\label{c1}	
	\end{align}
 
Our goal is to show that the limit $(\vr, \vt, \chi, \vu)$ is a weak solution of the same problem in the sense of Definition \ref{DMR1}. As expected, the proof shares some similarity with its counterpart for the Navier--Stokes--Fourier system discussed in detail in \cite[Chapters 2,3]{FeNo6A}. We therefore focus on the necessary modifications to accommodate the $\chi-$dependent transport coefficients.

\subsection{Compactness of the order parameter}

It follows from the strong convergence established in \eqref{c1}, the uniform bound \eqref{ap7}, and a simple interpolation argument that 
\begin{equation} \label{c2}
	\chi_n \to \chi \ \mbox{in}\ L^2(0,T; W^{1,2}(\Omega)). 
	\end{equation}

Consequently, it is a routine matter to perform the limit in the Allen--Cahn equation \eqref{MR14} obtaining 
	\begin{align} 
	\partial_t \chi + \vu \cdot \Grad \chi &= \Del \chi - \frac{\partial f}{\partial \chi } 
	\quad \mbox{a.a. in}\ (0,T) \times \Omega,\br \Grad \chi \cdot \vc{n}|_{\partial \Omega} &= 0,\  \chi(0, \cdot) = \chi_0.
	\label{c3}
\end{align}	

\subsection{Pointwise convergence of the temperature}

Our goal is to establish strong convergence of the temperatures $\vtn$. To see this, we first claim the uniform bound 
\begin{equation} \label{c4}
	\int_0^T  \| \log(\vt_n) \|_{W^{1,2}(\Omega)}^2 \dt 
	\end{equation}
that can be shown exactly as in \cite[Chapter 2, Section 2.2.4]{FeNo6A}. In particular, the bound \eqref{c4}  implies positivity of the absolute temperature a.a. in $(0,T) \times \Omega$.

In view of the bounds \eqref{ap4}, we have
\[
\vrn s(\vrn, \vtn) \ \mbox{bounded in}\ L^\infty(0,T; L^1(\Omega)) , 
\]
and 
\[
\vrn s(\vrn, \vtn) \ \mbox{bounded in}\ L^2(0,T; L^r(\Omega)) 
\ \mbox{for some}\ r > 1
\]
uniformly for $n \to \infty$.
Consequently, passing to a suitable subsequence,
\[
	\vrn s(\vrn, \vtn) \to \Ov{\vr s(\vr, \vt)} \ \mbox{weakly in}\ L^p((0,T) \times 
	\Omega) \ \mbox{for some}\ p > 1.
\] 
Here and hereafter, the bar denotes a weak limit of a composition of a weakly converging sequence with a (nonlinear) function.

The next step is using the structural hypotheses imposed on $s$ to deduce 
\[
|\vrn s(\vrn, \vtn) \vun | \aleq \Big( |\vun| |\vtn|^3 +\vrn |\log(\vrn)| |\vun| + |\vun| + \vrn |\log (\vtn)| |\vun|   \Big),
\]
where the uniform bounds \eqref{ap4} give rise to
\[
\Big( |\vun| |\vtn|^3 +\vrn |\log(\vrn)| |\vun| + |\vun| +  \vrn|\log (\vtn)| |\vun|  
\Big) \ \mbox{bounded in}\ L^r ((0,T) \times \Omega)
\]
for some $r > 1$. Here the most difficult term can be handled as 
\[
\|  \vrn  \log (\vtn) \vun \|_{L^{\frac{30}{29}}(\Omega; R^d)} \leq \| \vrn \vun \|_{L^{\frac{5}{4}}(\Omega; R^d)} \| \log (\vtn) \|_{L^6(\Omega)};
\]
whence 	\eqref{ap4}, together with \eqref{c2} and the Sobolev embedding 
$W^{1,2} \hookrightarrow L^6 (d=2,3)$  yield the desired conclusion. Thus we have 
\begin{equation} \label{c5}
	\vrn s(\vrn, \vtn) \vun \to \Ov{ \vr s(\vr, \vt ) \vu } \ \mbox{weakly in}\ 
	L^r((0,T) \times \Omega; R^d) \ \mbox{for some}\ r > 1.
\end{equation}

Writing the entropy flux in the form
\[
\frac{\kappa (\vtn, \chi_n) }{\vtn} |\Grad \vtn | \aleq \Big( |\Grad \log (\vtn) | + \vtn^{\frac{3}{2}} \left| \Grad \vtn^{\frac{3}{2}}
\right| \Big)
\]	
we may use the bounds  \eqref{ap4} to obtain 
\begin{equation} \label{c6}
	\frac{\kappa (\vtn, \chi_n) }{\vtn} \Grad \vtn \ \mbox{bounded in}\ L^r((0,T) \times \Omega; R^d) 
	\ \mbox{for some}\ r > 1
\end{equation}
uniformly for $n \to \infty$.

Now, going back to the entropy inequality \eqref{MR15}, we have
\begin{align}
	- \int_0^T \intO{ \vrn s(\vrn, \vtn) \partial_t \varphi } \dt
	&- \int_0^T \intO{ \vrn s (\vrn, \vtn) \vun \cdot \Grad \varphi } \dt
	\br &-
	\int_0^T \intO{ \frac{\vc{q}(\vt_n, \chin, \Grad \vtn) }{\vtn} \cdot \Grad \varphi } \dt  \geq 0
	\nonumber	
\end{align}
for any $\varphi \in C^1_c((0,T) \times \Ov{\Omega})$, $\varphi \geq 0$.

At this stage, we need the following version of Aubin--Lions Lemma, 
proved in 
\cite[Chapter 6, Lemma 6.3]{EF70}.

\begin{mdframed}[style=MyFrame]
	
	\begin{Lemma}[Lions--Aubin Lemma] \label{LionsAubin}
		
		\medskip
		
		Let $\{ v_n \}_{n=1}^\infty$ be a sequence of functions,
		\[
		\{ v_n \}_{n=1}^\infty \ \mbox{bounded in} \ L^2(0,T; L^q(\Omega)) \cap L^\infty(0,T; L^1(\Omega)),\ q > \frac{2 d}{d + 2}.
		\]
		In addition, suppose 
		\[
		\partial_t v_n \geq g_n \ \mbox{in}\ \mathcal{D}'((0,T) \times \Omega),
		\]  
		where

		\[
		\{ g_n \}_{n=1}^\infty \ \mbox{is bounded in}\ L^1(0,T ; W^{-m,r}(\Omega)) 
		\]
		for some $m \geq 1$, $r > 1$. 
		
		\medskip
		
		Then, up to a suitable subsequence, 
		\[
		v_n \to v \ \mbox{in}\ L^2(0,T; W^{-1,2}(\Omega)).
		\]

	\end{Lemma}
	
\end{mdframed}

Consequently, we get
\[
	\Ov{\vr s(\vr, \vt) g(\vt) } = \Ov{ \vr s(\vr, \vt) } \ \Ov{g(\vt)} 
\]
for any $g \in C^1_c (R)$. As a matter of fact, the above relation can be extended to any continuous non--decreasing function $g$ as soon as the quantities remain integrable.
In particular, 
\begin{equation} \label{c7}
	\Ov{\vr s(\vr, \vt) \vt } = \Ov{ \vr s(\vr, \vt) } \ \vt.
\end{equation}

Next, we use the renormalized equation \eqref{renorm} and apply the same argument to obtain 
\begin{equation} \label{c8}
	\Ov{b(\vr) g(\vt) } = \Ov{b(\vr)}\ \Ov{g(\vt)} 
\end{equation} 
for any bounded $b$.

Now, we recall \eqref{MR2} to write 
\[
\vrn s(\vrn, \vtn) = \vrn \mathcal{S} \left(  \frac{\vrn }{\vtn^{\frac{3}{2} }} \right) +  \frac{4a}{3} \vtn^3. 
\]
As $\mathcal{S}$ is decreasing (cf. \eqref{MR3}),  there holds
\[
\left[ \mathcal{S} \left(  \frac{\vrn }{\vtn^{\frac{3}{2} } } \right) -  \mathcal{S} \left(  \frac{\vrn }{\vt^{\frac{3}{2} } } \right) \right] (\vtn - \vt ) \geq 0.
\] 
Further, we claim that 
\[
\int_0^T \intO{ \vr_n \mathcal{S} \left(  \frac{\vrn }{\vt^{\frac{3}{2} } } \right)  (\vtn - \vt ) } \dt \to 0. 
\]
Indeed this can be shown exactly as in \cite[Chapter 3, Section 3.7.3]{FeNo6A}, using the renormalized equation \eqref{renorm} and boundedness of the temperature gradients. 

Consequently, we get 
\[
\Ov{ \vr \mathcal{S} \left(  \frac{\vr }{\vt^{\frac{3}{2} } } \right) \vt } \geq 
\Ov{ \vr \mathcal{S} \left(  \frac{\vr }{\vt^{\frac{3}{2} } } \right)} \vt ,
\]
and, similarly,  
\[
	\Ov{ \vt^3 \vt } \geq \Ov{\vt^3} \vt.
\]
Thus going back to \eqref{c7} we conclude 
\[
\Ov{ \vt^3 \vt } = \Ov{\vt^3} \vt,
\]
which yields, modulo a suitable subsequence, 
\begin{equation} \label{c10}
	\vtn \to \vt \ \mbox{a.a. in}\ (0,T) \times \Omega.
	\end{equation}

\subsection{Pointwise convergence of the densitites}
\label{pe}

Our ultimate goal is to establish strong (pointwise) convergence of the densities.

\subsubsection{Pressure estimates} 

We start by testing the momentum 
equation \eqref{MR13} on 
\[
\bfphi (t,x) = \phi \Grad \Del^{-1} [\phi  b(\vr)],\  \ \phi \in \DC(\Omega),
\]  
where $\Del^{-1}$ is the inverse of the Laplacean defined on $R^d$ via the convolution with the Poisson kernel.
The resulting expression reads
\begin{align} 
	\int_0^{T} &\intO{ p(\vrn, \vtn)  \left[ \phi^2 b(\vrn) + \Grad \phi \cdot \Grad \Del^{-1} [\phi  b(\vrn)] \right]  } \dt \br &=
	\left[ \intO{ \phi \vrn \vun \cdot \Grad \Del^{-1} [\phi  b(\vrn)] }  \right]_{t=0}^{t = T} \br
	&- \int_0^{T} \intO{ \phi \vrn \vun \cdot \partial_t \Big( \Grad \Del^{-1} [\phi  b(\vrn)] \Big) } \dt\br
	&- \int_0^{T} \intO{ \phi \vrn \vun \otimes \vun : \Grad^2 \Del^{-1} [\phi  b(\vrn)] } \dt\br
	&- \int_0^{T} \intO{ \vrn \vun \otimes \vun \cdot \Grad \phi \cdot \Grad \Del^{-1} [\phi  b(\vrn)] } \dt
	\br
	&+ \int_0^{T} \intO{ \phi \mathbb{S}  : \Grad^2 \Del^{-1} [\phi  b(\vrn)] } \dt \br
	&- \int_0^T \intO{ \left( \Grad \chi_n \otimes \Grad \chi_n - \frac{1}{2} |\Grad \chi_n|^2 \mathbb{I} \right) : \Grad^2 \Del^{-1} [\phi  b(\vrn)] } \dt \br 
	&+\int_0^{T} \intO{ \mathbb{S} \cdot \Grad \phi \cdot \Grad \Del^{-1} [\phi  b(\vrn)] } \dt \br 
	&- \int_0^T \intO{ \left( \Grad \chi_n \otimes \Grad \chi_n - \frac{1}{2} |\Grad \chi_n|^2 \mathbb{I} \right) \cdot \Grad \phi \cdot \Grad \Del^{-1} [\phi  b(\vrn)]} \dt \br 
	&+ \int_0^{T} \intO{ f(\chin)  \left[ \phi^2 b(\vrn) + \Grad \phi \cdot \Grad \Del^{-1} [\phi  b(\vrn)] \right]  } \dt
	\label{c11}
\end{align}
where, by virtue of the renormalized equation of continuity \eqref{renorm},
\begin{equation} \label{c12}
	\begin{split}
		\partial_t &   \Big( \Grad \Del^{-1} [\phi  b(\vrn)] \Big) \\
		&= - \Grad \Del^{-1} \Big[ \phi \Div ( b(\vrn) \vun) \Big] 
		+ \Grad \Del^{-1} \Big[ \phi (b(\vrn) - b'(\vrn) \vrn) \Div \vun \Big]. 
	\end{split}
\end{equation}

In view of the energy estimates \eqref{ap4}, the integrals on the right--hand side are already bounded as long as 
$b(\vr) \approx \vr^\beta$, with $\beta > 0$ small enough. In particular, the relation \eqref{c11} implies the desired pressure estimates, 
\begin{equation} \label{c13}
	\int_0^T \int_K p(\vrn, \vtn) \vrn^\beta \dx \ \dt \leq c(K) \ \mbox{for any compact}\ K \subset \Omega.
	\end{equation}

The local estimate \eqref{c13} can be extended up to the boundary by considering the test function 
\[
\bfphi (t,x) = 
\mathcal{B} \left[ \Phi  \right],\ \Phi \in L^q(\Omega), \ \intO{ \Phi } = 0,
\]
where $\mathcal{B}$ is the so--called Bogovskii operator, see \cite[Appendix, Section 11.6]{FeNo6A}. Similarly to the above, we obtain
\begin{align} 
		\int_0^{T} &\intO{ \Phi p (\vrn, \vtn)  } \dt =
		\left[ \intO{ \vrn \vun \cdot \mathcal{B} \left[ \Phi \right] }  \right]_{t=0}^{t = T} \br
		&- \int_0^{T} \intO{  \vrn \vun \otimes \vun : \Grad \mathcal{B} \left[ \Phi \right] } \dt
		+ \int_0^{T} \intO{  \mathbb{S}  : \Grad \mathcal{B} \left[\Phi \right] } \dt \br
		&-
 \int_0^{T} \intO{ \left( \Grad \chin \otimes \Grad \chin - \frac{1}{2} |\Grad \chi_n |^2 \mathbb{I} \right)  : \Grad \mathcal{B} \left[\Phi \right] } \dt \br	
		&- \int_0^{T} \intO{  f(\chin)  \Phi  } \dt.
	\label{c14}
\end{align}

A suitable choice of $\Phi$ gives rise to the estimate 
\begin{equation} \label{c15}
	\int_0^T \intO{ p(\vrn, \vtn) {\rm dist}^{- \omega } [x, \partial \Omega] } \dt \aleq 1 \ \mbox{for some}\ \omega > 0,
\end{equation}

The estimates \eqref{c13}, \eqref{c15} imply equi--integrability of the pressure in $(0,T) \times \Omega$. Alternatively, one can use directly the Bogovksii operator to obtain integrability of the pressure in the Lebesgue space $L^r$ with $r > 1$, 
	see \cite[Chapter 3, Section 3.7]{FeNo6A} for details. 
	
	\subsubsection{Convergence}
	
	With pressure estimates at hand, we can perform the limit in the momentum equation \eqref{M13}: 
\begin{align} 
	\int_0^T \intO{ \Big[ &\vr \vu \cdot \partial_t \bfphi + \vr \vu \otimes \vu : \Grad \bfphi + \Ov{p(
			\vr, \vt) } \Div \bfphi \Big] } \dt = \int_0^T \intO{ \mathbb{S} (\vt, \chi, \Ds \vu) : \Grad \bfphi } \dt \br 
	&- \int_0^T \intO{ \left[ \Grad \chi \otimes \Grad \chi - \frac{1}{2} |\Grad \chi |^2 \mathbb{I} \right] : \Grad \bfphi } \dt  + \int_0^T \intO{ f(\chi) \Div \bfphi } \dt \br &= - 
	\intO{ (\vr \vu)_0 \cdot \bfphi (0, \cdot) }
	\label{c16}
\end{align}
for any $\bfphi \in C^1_c([0,T) \times \Omega; R^d)$. Note that we have already established pointwise convergence of $\{ 
\vt_n, \chi_n, \Grad \chi_n \}_{n  \geq 0}$.

Let 
\[
T_k(r) = k T \left( \frac{r}{k} \right),\ r \geq 0,\ k \geq 1,
\]
\[
T \in C^\infty[0,\infty),\ 
T(z) = \left\{ \begin{array}{l} z \ \mbox{if}\ 0 \leq z \leq 1, \\ 
	\ \mbox{concave on}\ [0, \infty), \\ 2 \ \mbox{if}\ f \geq 3 
\end{array}  \right.
\]
be a cut--off function. 

Following the arguments of \cite[Chapter 3, Section 3.7.4]{FeNo6A} we first use the identity \eqref{c11} with 
$b = T_k$ and let $n \to \infty$. Then we consider 
\[
\bfphi (t,x) = \phi \Grad \Del^{-1} [\phi \Ov{T_k(\vr)}] 
\]
as a test function in the limit momentum equation \eqref{c16}. Comparing the limits and performing the arguments of 
\cite[Chapter 3, Section 3.7.4]{FeNo6A} we arrive at the identity 
\begin{align} 
		\int_0^T &\intO{ \phi^2 \left( \Ov{ p (\vr, \vt) T_k (\vr)} - \Ov{ p (\vr, \vt)}\ \Ov{T_k(\vr)} \right) } \dt \br 
		&= \lim_{n \to \infty} \int_0^T \int_\Omega \phi \Big( \mathbb{S}(\vtn, \chin, \Ds \vun) : \Grad^2 \Del^{-1}[\phi T_k (\vrn)] - \mathbb{S}(\vt, \chi,  \Ds \vu) : \Grad^2 \Del^{-1}[\phi \Ov{T_k (\vr)}] 
		\Big) \dx \dt.            
	\label{c17}
\end{align}

The rightmost integral in \eqref{c17} can be rewritten as
\begin{align}
	\int_\Omega &\phi \Big( \mathbb{S}(\vtn, \chin, \Ds \vun) : \Grad^2 \Del^{-1}[\phi T_k (\vrn)] \br &- \mathbb{S}(\vt, \chi, \Ds \vu) : \Grad^2 \Del^{-1}[\phi \Ov{T_k (\vr)}] 
	\Big) \dx  \br
	&= \int_\Omega \phi \Big( \Grad^2 \Del^{-1} : \Big[ \phi \mathbb{S}(\vtn, \chin, \Ds \vun) \Big] T_k (\vrn) \br &- \Grad^2 \Del^{-1}: \Big[ \phi \mathbb{S}(\vt, \chi, \Ds \vu) \Big] \Ov{T_k (\vr)} 
	\Big) \dx  
	\nonumber
\end{align}
Next,  
\[
\begin{split}
	\Grad^2 \Del^{-1} &: \Big[ \phi \mathbb{S}(\vtn, \chin, \Ds \vun) \Big] = 
	\phi \left( \frac{2}{3} \mu(\vtn, \chi_n) + \eta(\vtn, \chin) \right) \Div \vun  \\
	&+  \Grad^2 \Del^{-1} : \Big[ \phi \mathbb{S}(\vtn, \chin, \Ds \vun) \Big] - 
	\phi \left( \frac{2}{3} \mu(\vtn, \chin) + \eta(\vtn, \chin) \right) \Div \vun,
\end{split}
\]
and, similarly, 
\[
\begin{split}
	\Grad^2 \Del^{-1} &: \Big[ \phi \mathbb{S}(\vt, \chi, \Ds \vu) \Big] = 
	\phi \left( \frac{2}{3} \mu(\vt, \chi) + \eta(\vt, \chi) \right) \Div \vu \\
	&+  \Grad^2 \Del^{-1} : \Big[ \phi \mathbb{S}(\vt, \chi, \Ds \vu) \Big] - 
	\phi \left( \frac{2}{3} \mu(\vt, \chi) + \eta(\vt, \chi) \right) \Div \vu.
\end{split}
\]

At this stage, we use 
the following result in the spirit of Coifman and Meyer \cite{COME}, see \cite[Chapter 11, Section 11.18, Theorem 11.35]{FeNo6A}.
Let
\[
\mathcal{R} = \mathcal{R}_{i,j = 1}^d,\ 
\mathcal{R}_{i,j} [v] (x) = \mathcal{F}^{-1}_{\xi \to x} \left[ \frac{\xi_i \xi_j }{|\xi|^2 } \mathcal{F}_{x \to \xi }[v] \right]
\]
where $\mathcal{F}: R^d \to R^d$ denotes the Fourier transform.

\begin{mdframed}[style=MyFrame]
	
	\begin{Lemma}[Commutator Lemma] \label{CommL}
		
		\medskip

		Let $w \in W^{1,r}(R^d)$, $v \in L^p(R^d; R^d)$, 
		\[
		1 < r < d,\ 1 < p < \infty,\ \frac{1}{r} + \frac{1}{p} - \frac{1}{d} < 1
		\]
		be given.

		Then for any $s$ satisfying 
		\[
		\frac{1}{r} + \frac{1}{p} - \frac{1}{d} < \frac{1}{s} < 1, 
		\]
		there exists

		\[
		0 < \beta < 1,\ \frac{\beta}{d} = \frac{1}{s} + \frac{1}{d} - \frac{1}{p} - \frac{1}{r} 
		\]
		such that 
		\[
		\Big\| \mathcal{R}[ w \vc{v} ] - w \mathcal{R} [\vc{V}] \Big\|_{W^{\beta,s}(R^d; R^d)} \leq 
		c \| w \|_{W^{1,r}(R^d)} \| \vc{V} \|_{L^p(R^d; R^d)}, 
		\]
		with a positive constant $c = c(s,p,r)$.
		
	\end{Lemma}
	
\end{mdframed}

Applying Lemma \ref{CommL}, we rewrite relation \eqref{c17} in the form 
\begin{equation} \label{c18}
	\begin{split}
		\int_0^T &\intO{ \phi^2 \left( \Ov{ p (\vr, \vt) T_k (\vr)} - \Ov{(p(\vr, \vt)}\ \Ov{T_k(\vr)} \right) } \dt \\ 
		&= \int_0^T  \intO{ \phi^2 \left( \frac{2}{3} \mu(\vt) + \eta(\vt) \right) \Big( \Ov{T_k(\vr) \Div \vu } - \Ov{T_k(\vr) }
			\Div \vu \Big)}  \dt          
	\end{split}	
\end{equation}
for any $\phi \in \DC(\Omega)$. 

Having established \eqref{c18}, the rest of the proof of strong convergence of the density is the same as in  \cite[Chapter 11, Section 11.18, Theorem 11.35]{FeNo6A}. Thus we get 
\begin{equation} \label{c19}
	\vr_n \to \vr \ \mbox{a.a. in}\ (0,T) \times \Omega 
	\end{equation}
as soon as 
\[
\vr_n (0, \cdot) \to \vr(0, \cdot) \ \mbox{in}\ L^{\frac{5}{3}}(\Omega).
\]

\section{Existence of weak solutions}
\label{E}

The existence of weak solutions can be shown by means of a suitable modification of the approximate scheme introduced in 
\cite[Chapter 3]{FeNo6A}. 

First consider a system of functions $\{ \vc{w}_j \}_{j=1}^\infty \subset \DC(\Omega; R^d)$ that forms an orthonormal basis of the space $L^2(\Omega; R^d)$ and fix two positive parameters $\ep > 0$, $\delta > 0$.
The approximate velocity is looked for in the form 
\[
\vu \in C^1([0,T]; X_n),\ X_n = {\rm span} \left\{ \vc{w}_j \ \Big| \ 1 \leq j \leq n \right\},
\] 
whereas the functions $(\vr, \vt, \chi, \vu)$ solve the following system of equations: 

\begin{itemize}
	\item {\bf Vanishing viscosity approximation of the equation of continuity.}
	\begin{align} 
		\partial_t \vr + \Div (\vr \vu) &= \ep \Del \vr \ \mbox{in}\ (0,T) \times \Omega, \br
		\ep \Grad \vr \cdot \vc{n} &= 0\ \mbox{in}\ (0,T) \times \partial \Omega, \br 
		\vr(0, \cdot) &= \vr_{0,\delta}.
		\label{E1}	
	\end{align}
	
	\item {\bf Galerkin approximation of the momentum equation.}	
	
	\begin{align}
		\int_0^\tau &\intO{ \Big[ \vr \vu \cdot \partial_t \bfphi + \vr \vu \otimes \vu : \Grad \bfphi + \delta \left( \vr^\Gamma + \vr^2 \right) \Div \bfphi + p(\vr, \vt) \Div \bfphi \Big] } \dt \br &= 
		\int_0^\tau \intO{ \Big[ \mathbb{S}_\delta (\vt, \Ds \vu) : \Ds \bfphi 
			- \left( \Grad \chi \otimes \Grad \chi - \frac{1}{2} |\Grad \chi |^2 \mathbb{I} \right) : \Ds \bfphi 
			+ f(\chi) \Div \bfphi  \Big] } \dt \br & + \left[ \intO{ \vr \vu \cdot \bfphi } \right]_{t=0}^{t = \tau}, \ 
		\vr \vu (0, \cdot) = \vm_0,  
		\label{E2}
	\end{align}
	for any $0 \leq \tau \leq T$, $\bfphi \in C^1([0,T]; X_n)$. Here $\mathbb{S}_\delta$ denotes the viscous stress with the shear viscosity $\mu_\delta = \mu + \delta \vt$.

	\item {\bf Approximate internal energy balance.}
	
	\begin{align}
		\partial_t \Big[ \vr (e(\vr, \vt) + \delta \vt) \Big] &+ \Div \Big[ \vr (e(\vr, \vt) + \delta  \vt) \vu \Big] - \Div \left[\left(\delta \left( \vt^\Gamma + \frac{1}{\vt} \right)  +    \kappa (\vt, \chi) \right) \Grad \vt \right] \br & = \mathbb{S}_\delta : \Ds \vu - 
		p(\vr, \vt) \Div \vu +  \left( \Del \chi -  \frac{\partial f(\chi) }{\partial \chi } \right)^2 \br &+ \ep \delta \left( \Gamma \vr^{\Gamma - 2} + 2 \right) |\Grad \vr|^2 + \frac{\delta}{\vt^2} - \ep \vt^5 \ \mbox{in}\ (0,T) \times \Omega, \br
		\Grad \vt \cdot \vc{n} &= 0 \ \mbox{in}\ (0,T) \times \partial \Omega, \br
		\vt(0, \cdot) &= \vt_0.
		\label{E3}
	\end{align}	

\item {\bf Allen--Cahn equation.} 
\begin{align}
	\partial_t \chi + \vu \cdot \Grad \chi &= \Del \chi - \frac{f(\chi) }{\partial \chi } \ \mbox{in}\ (0,T) \times \Omega, \br
	\Grad \chi \cdot \vc{n} &= 0 \ \mbox{in}\ (0,T) \times \partial \Omega, \br
	\chi(0, \cdot) &= \chi_0.
	\label{E4}
	\end{align}
	
\end{itemize}

Given the available estimates established in Section \ref{ap}, the compactness arguments of Section \ref{c} can be adapted to prove Theorem \ref{TMR1} performing successively the limits $n \to 0$, $\ep \to 0$, and $\delta \to 0$.

\bigskip 

\centerline{\bf Acknowledgement}

\bigskip

The work was supported by the Grant 8J20FR007 in the framework of PHC Barrande of France and Czech Republic.

\def\cprime{$'$} \def\ocirc#1{\ifmmode\setbox0=\hbox{$#1$}\dimen0=\ht0
	\advance\dimen0 by1pt\rlap{\hbox to\wd0{\hss\raise\dimen0
			\hbox{\hskip.2em$\scriptscriptstyle\circ$}\hss}}#1\else {\accent"17 #1}\fi}


\end{document}